\documentclass[11pt]{amsart}

\usepackage[margin=1in]{geometry}
\usepackage{amsmath,amssymb,amsthm}
\usepackage{mathptmx}
\usepackage[T1]{fontenc}
\usepackage{enumitem}
\usepackage[colorlinks=true,linkcolor=blue,citecolor=blue,urlcolor=blue]{hyperref}

\newtheorem{theorem}{Theorem}[section]
\newtheorem{proposition}[theorem]{Proposition}
\newtheorem{lemma}[theorem]{Lemma}
\newtheorem{corollary}[theorem]{Corollary}
\theoremstyle{definition}
\newtheorem{definition}[theorem]{Definition}
\newtheorem{example}[theorem]{Example}
\theoremstyle{remark}
\newtheorem{remark}[theorem]{Remark}

\newcommand{\R}{\mathbb{R}}
\newcommand{\Z}{\mathbb{Z}}
\newcommand{\cK}{\mathcal{K}}
\DeclareMathOperator{\spn}{span}

\title[No nontrivial ACR in weakly reversible networks]{Weakly reversible mass-action systems have no\\ nontrivial absolute concentration robustness}

\author{Sibo Peng}
\address{Huazhong University of Science and Technology, Wuhan 430074, Hubei, China}
\email{kuhaizhongdebailanren@gmail.com}

\author{Chen Deng}
\address{Huazhong University of Science and Technology, Wuhan 430074, Hubei, China}
\email{sunnyorange.andrew@gmail.com}

\author{Yuntao Zou}
\address{Huazhong University of Science and Technology, Wuhan 430074, Hubei, China}
\email{zouyuntao@hust.edu.cn}

\author{Ruixuan Li}
\address{Huazhong University of Science and Technology, Wuhan 430074, Hubei, China}
\email{rxli@hust.edu.cn}

\date{\today}


\begin{document}

\begin{abstract}
A mass-action system has absolute concentration robustness (ACR) in a species if that species takes the same value at every positive steady state. The robustness is trivial when the species appears in no conservation law, and nontrivial otherwise. All known networks with nontrivial ACR, including those produced by the criterion of Shinar and Feinberg, fail to be weakly reversible. We prove that this is necessary. If a weakly reversible network with fixed positive rate constants has a species $X_i$ that appears in some conservation law, then the value of $x_i$ over the positive steady states has infimum $0$ and supremum $+\infty$. In particular, weakly reversible mass-action systems admit only trivial ACR. The proof combines two ingredients. The first is an a priori bound: within a fixed stoichiometric compatibility class, the positive steady states stay in a compact box, uniformly over all independent rescalings of the rate constants of the individual linkage classes. The second is the classical correspondence between such rescalings and translations of the steady-state set along the orthogonal complement of the stoichiometric subspace. The only external result used is the existence theorem of Boros for weakly reversible systems. The main results have been formally verified in Lean~4, taking the existence theorem of Boros as an explicit hypothesis.
\end{abstract}

\maketitle

\section{Introduction}

Consider a mass-action reaction network with species $X_1,\dots,X_n$ and fixed positive rate constants. The system has \emph{absolute concentration robustness} (ACR) in $X_i$ if it has at least one positive steady state and the concentration $x_i$ is the same at all of them. Shinar and Feinberg introduced the notion and gave a structural sufficient condition: a deficiency-one network with a positive steady state, containing two nonterminal complexes that differ only in $X_i$, has ACR in $X_i$ \cite{SF10}. Biological examples include the EnvZ--OmpR osmoregulatory system \cite{SF10,BG03}. Since then ACR has been studied from algebraic, geometric and dynamical points of view \cite{JC22,MST22,KNS24}.

If $X_i$ appears in no conservation law, that is, if the unit vector $e_i$ lies in the stoichiometric subspace $S$, then nothing ties $x_i$ to a conserved total and ACR in $X_i$ carries little structural information. We call ACR in $X_i$ \emph{trivial} if $e_i\in S$ and \emph{nontrivial} otherwise. The Shinar--Feinberg networks have nontrivial ACR, and they always contain nonterminal complexes, so they are never weakly reversible. This raises the question of whether weakly reversible networks can have nontrivial ACR at all. The answer is known to be negative for complex-balanced, in particular deficiency-zero, weakly reversible systems, by the toric structure of complex-balanced equilibria \cite{HJ72} (see also \cite{JC22}), and the same argument applies whenever the linkage classes are stoichiometrically independent. The difficulty lies in networks whose linkage classes are stoichiometrically dependent, where steady states need not be complex balanced and the steady-state set has no evident toric structure.

Our main result settles the question for all weakly reversible networks and all rate constants.

\begin{theorem}\label{thm:main}
Let $(G,\kappa)$ be a weakly reversible mass-action system with $\kappa>0$, and let $E_\kappa$ be its set of positive steady states. If $e_i\notin S$, then
\[
\inf_{x\in E_\kappa} x_i = 0 \qquad\text{and}\qquad \sup_{x\in E_\kappa} x_i = +\infty .
\]
\end{theorem}

\begin{corollary}\label{cor:acr}
A weakly reversible mass-action system has ACR in $X_i$ only if $e_i\in S$. In other words, weakly reversible mass-action systems admit no nontrivial ACR.
\end{corollary}

Theorem~\ref{thm:main} holds for every positive rate vector, not only for generic ones, and it gives more than the failure of ACR: across the steady states, $x_i$ comes arbitrarily close to $0$ and becomes arbitrarily large. Weak reversibility cannot be dropped, as the example in Section~\ref{sec:examples} shows.

The proof has two steps. In Section~\ref{sec:bound} we show that, in a fixed compatibility class, the positive steady states of a weakly reversible system lie in a compact box $[m,M]^n$ that does not depend on how the rate constants of the different linkage classes are rescaled relative to each other (Proposition~\ref{prop:uniform}). The argument looks at the lowest level of a linear functional on each linkage class; the reaction that leaves this level upward dominates every reaction that comes down, whatever the rescaling. In Section~\ref{sec:proof} we combine this bound with the well-known identity that rescaling the rate constants of linkage class $K$ by $e^{t\psi_K}$, for $\psi\in S^\perp$, is equivalent to translating steady states by $e^{t\psi}$. The existence theorem of Boros \cite{Bor19} provides a steady state in the box for each rescaling, and translating it sweeps $x_i$ across $(0,\infty)$.

\section{Setting}\label{sec:setting}

A reaction network $G$ consists of finitely many reactions $r\colon y_r\to y'_r$, where the complexes $y_r\neq y'_r$ lie in $\Z^n_{\ge 0}$. Write $\Delta_r=y'_r-y_r$ and
\[
S=\spn\{\Delta_r\}\subseteq \R^n .
\]
The \emph{complex graph} has the complexes as vertices and the reactions as directed edges. Its connected components (ignoring orientation) are the \emph{linkage classes}, which we denote by $K\in\cK$. The network is \emph{weakly reversible} if every linkage class is strongly connected.

Given rate constants $\kappa=(\kappa_r)>0$, the mass-action vector field is
\[
f_\kappa(x)=\sum_r \kappa_r x^{y_r}\Delta_r,\qquad x^{y}=\prod_{k=1}^n x_k^{y_k},
\]
and the set of positive steady states is $E_\kappa=\{x\in\R^n_{>0}: f_\kappa(x)=0\}$. For $b\in\R^n_{>0}$ the positive stoichiometric compatibility class of $b$ is $P_b=(b+S)\cap\R^n_{>0}$. We split $f_\kappa$ by linkage class,
\[
f_\kappa=\sum_{K\in\cK} f_K,\qquad f_K(x)=\sum_{r\in K}\kappa_r x^{y_r}\Delta_r,
\]
where $r\in K$ means that the reaction $r$ belongs to the linkage class $K$. For $\alpha=(\alpha_K)\in\R^{\cK}_{>0}$ define the rescaled field
\[
F_\alpha=\sum_{K\in\cK}\alpha_K f_K ,
\]
which is the mass-action field of the same network with rate constants $\alpha_K\kappa_r$ for $r\in K$.

\begin{definition}
The system $(G,\kappa)$ has \emph{absolute concentration robustness} (ACR) in $X_i$ if $E_\kappa\neq\emptyset$ and $x_i$ is constant on $E_\kappa$. The ACR is \emph{trivial} if $e_i\in S$ and \emph{nontrivial} otherwise.
\end{definition}

Since $S^\perp$ is the orthogonal complement of $S$, the condition $e_i\notin S$ is equivalent to the existence of $\psi\in S^\perp$ with $\psi_i\neq 0$, that is, of a conservation law involving $X_i$.

We use the following existence theorem.

\begin{theorem}[Boros \cite{Bor19}, Theorem 1]\label{thm:boros}
If $G$ is weakly reversible and $\kappa>0$, then every positive stoichiometric compatibility class contains at least one positive steady state.
\end{theorem}

The next identity is standard; we include its one-line proof.

\begin{lemma}\label{lem:shift}
Let $\psi\in S^\perp$. Then $\psi\cdot y$ is constant on each linkage class $K$; denote this constant by $\psi_K$. For every $x\in\R^n_{>0}$ and $t\in\R$,
\[
f_\kappa\big(e^{t\psi}\circ x\big)=F_{\alpha(t)}(x),\qquad \alpha_K(t)=e^{t\psi_K},
\]
where $e^{t\psi}\circ x=(e^{t\psi_1}x_1,\dots,e^{t\psi_n}x_n)$.
\end{lemma}

\begin{proof}
Each reaction joins two complexes of the same linkage class, and $\psi\cdot\Delta_r=0$ for all $r$, so $\psi\cdot y$ is constant along every edge and hence on every linkage class. For $r\in K$ we have $(e^{t\psi}\circ x)^{y_r}=e^{t\psi\cdot y_r}x^{y_r}=e^{t\psi_K}x^{y_r}$, and summing over $r$ gives the claim.
\end{proof}

\section{A bound uniform in the linkage-class rescalings}\label{sec:bound}

\begin{proposition}\label{prop:uniform}
Let $G$ be weakly reversible, $\kappa>0$ and $b\in\R^n_{>0}$. There exist $0<m<M<\infty$ such that, for every $\alpha\in\R^{\cK}_{>0}$, every positive zero of $F_\alpha$ in $P_b$ lies in $[m,M]^n$.
\end{proposition}

\begin{proof}
Suppose not. Then there are $\alpha^{(j)}\in\R^{\cK}_{>0}$ and $x^{(j)}\in P_b$ with $F_{\alpha^{(j)}}(x^{(j)})=0$ and
\[
L_j:=\|\log x^{(j)}\|_\infty\longrightarrow\infty ,
\]
where $\log$ acts coordinatewise. Put $u^{(j)}=-\log x^{(j)}/L_j$, so that $\|u^{(j)}\|_\infty=1$ and $x^{(j)}_k=\exp(-L_j u^{(j)}_k)$. Passing to a subsequence we may assume $u^{(j)}\to u$ with $\|u\|_\infty=1$.

\emph{Step 1: the sign of $u\cdot f_K$.} Fix a linkage class $K$. If $u\cdot y$ is constant on the complexes of $K$, then $u\cdot\Delta_r=0$ for every $r\in K$, and $u\cdot f_K\equiv 0$.

Otherwise let $h_0$ be the minimum of $u\cdot y$ over the complexes of $K$. Since $u\cdot y$ is not constant on $K$ and $K$ is strongly connected, there is a directed path from a complex at level $h_0$ to a complex above level $h_0$. The first reaction $r_*$ on this path that leaves level $h_0$ has source $y_{r_*}$ at level $h_0$ and satisfies $u\cdot\Delta_{r_*}>0$. If $r\in K$ satisfies $u\cdot\Delta_r<0$, then $u\cdot y_r>u\cdot y'_r\ge h_0$, so $u\cdot(y_r-y_{r_*})>0$. For such $r$,
\[
\frac{(x^{(j)})^{y_r}}{(x^{(j)})^{y_{r_*}}}=\exp\!\big(-L_j\,u^{(j)}\cdot(y_r-y_{r_*})\big)\longrightarrow 0,
\]
because $u^{(j)}\cdot(y_r-y_{r_*})\to u\cdot(y_r-y_{r_*})>0$ and $L_j\to\infty$. Dividing by $(x^{(j)})^{y_{r_*}}$,
\[
\frac{u\cdot f_K(x^{(j)})}{(x^{(j)})^{y_{r_*}}}
\;\ge\; \kappa_{r_*}\,u\cdot\Delta_{r_*}\;+\sum_{r\in K:\,u\cdot\Delta_r<0}\kappa_r\,(u\cdot\Delta_r)\,\frac{(x^{(j)})^{y_r}}{(x^{(j)})^{y_{r_*}}},
\]
since the reactions with $u\cdot\Delta_r\ge 0$ other than $r_*$ contribute nonnegative terms. The first term on the right is a fixed positive number and the sum tends to $0$. Hence $u\cdot f_K(x^{(j)})>0$ for all large $j$.

\emph{Step 2: $u\in S^\perp$.} There are finitely many linkage classes, so for all large $j$ every class satisfies either $u\cdot f_K\equiv0$ or $u\cdot f_K(x^{(j)})>0$. Since
\[
0=u\cdot F_{\alpha^{(j)}}(x^{(j)})=\sum_{K}\alpha^{(j)}_K\,u\cdot f_K(x^{(j)})
\]
and every $\alpha^{(j)}_K$ is positive, no class can be of the second kind. Thus $u\cdot y$ is constant on every linkage class, so $u\cdot\Delta_r=0$ for all $r$, that is, $u\in S^\perp$.

\emph{Step 3: contradiction.} As $x^{(j)}-b\in S$ and $u\in S^\perp$, we have $u\cdot(x^{(j)}-b)=0$ for all $j$. On the other hand, if $u_k>0$ then $x^{(j)}_k=\exp(-L_ju^{(j)}_k)\to0$, so $u_k(x^{(j)}_k-b_k)\to -u_kb_k<0$. If $u_k<0$ then $x^{(j)}_k\to\infty$ and $u_k(x^{(j)}_k-b_k)<0$ for large $j$. Terms with $u_k=0$ vanish. Since $u\neq0$, we get $u\cdot(x^{(j)}-b)<0$ for all large $j$, a contradiction.
\end{proof}

\begin{remark}\label{rem:compact}
Taking $\alpha\equiv 1$, Proposition~\ref{prop:uniform} shows that the positive steady states of a weakly reversible system in a fixed compatibility class form a compact subset of $\R^n_{>0}$. We do not claim novelty for this special case. The point needed below is the uniformity in $\alpha$, which comes for free from the argument: the sign of each $u\cdot f_K(x^{(j)})$ is independent of the positive factor $\alpha_K$.
\end{remark}

\begin{remark}\label{rem:weaker}
Weak reversibility enters the proof of Proposition~\ref{prop:uniform} only through the following property of each linkage class $K$: for every $u$ such that $u\cdot y$ is not constant on $K$, some reaction of $K$ has its source at the minimum of $u\cdot y$ over the complexes of $K$ and satisfies $u\cdot\Delta_r>0$. Theorem~\ref{thm:main} therefore holds for any network with this property for which every rescaled system $F_\alpha$ has a positive steady state in every compatibility class. We do not pursue this generalization.
\end{remark}

\section{Proof of the main theorem}\label{sec:proof}

\begin{proof}[Proof of Theorem~\ref{thm:main}]
Assume $e_i\notin S$. Let $\psi$ be the orthogonal projection of $e_i$ onto $S^\perp$. Then $\psi\in S^\perp$ and
\[
\psi_i=\psi\cdot e_i=\|\psi\|^2>0 .
\]
Let $\psi_K$ be the constant value of $\psi\cdot y$ on the linkage class $K$ (Lemma~\ref{lem:shift}). Fix $b\in\R^n_{>0}$ and let $m,M$ be as in Proposition~\ref{prop:uniform}.

For $t\in\R$ let $\alpha_K(t)=e^{t\psi_K}$. The system $F_{\alpha(t)}$ is the mass-action system of the same weakly reversible network with positive rate constants $\alpha_K(t)\kappa_r$. By Theorem~\ref{thm:boros} it has a positive steady state $\xi_t\in P_b$, and by Proposition~\ref{prop:uniform}, $\xi_t\in[m,M]^n$.

Set $z_t=e^{t\psi}\circ\xi_t$. By Lemma~\ref{lem:shift},
\[
f_\kappa(z_t)=F_{\alpha(t)}(\xi_t)=0,
\]
so $z_t\in E_\kappa$. Its $i$-th coordinate satisfies
\[
m\,e^{t\psi_i}\le (z_t)_i=e^{t\psi_i}(\xi_t)_i\le M\,e^{t\psi_i}.
\]
Since $\psi_i>0$, letting $t\to-\infty$ gives $\inf_{E_\kappa}x_i=0$, and letting $t\to+\infty$ gives $\sup_{E_\kappa}x_i=+\infty$.
\end{proof}

\begin{proof}[Proof of Corollary~\ref{cor:acr}]
If $e_i\notin S$, Theorem~\ref{thm:main} shows that $x_i$ is not constant on $E_\kappa$, so there is no ACR in $X_i$.
\end{proof}

\begin{remark}
The proof gives slightly more. For every $\psi\in S^\perp$ and every compatibility class $P_b$, the steady-state set $E_\kappa$ contains a curve $t\mapsto z_t$ (not necessarily continuous) with $\log z_t-t\psi$ bounded. Thus $E_\kappa$ meets every direction of $S^\perp$ at logarithmic scale, which is the asymptotic form of the toric structure that complex-balanced systems have exactly.
\end{remark}

\section{Examples}\label{sec:examples}

\begin{example}[Weak reversibility cannot be dropped]
Consider the network $B\to A$, $A+B\to 2B$ with rate constants $\kappa_1,\kappa_2$. Its stoichiometric subspace is spanned by $(1,-1)$, so $A$ appears in the conservation law $a+b=\text{const}$. The steady-state equation is $b(\kappa_1-\kappa_2a)=0$, so every positive steady state has $a=\kappa_1/\kappa_2$, and the system has nontrivial ACR in $A$. The network has two linkage classes, $\{B,A\}$ and $\{A+B,2B\}$, neither strongly connected. This is the simplest Shinar--Feinberg network \cite{SF10}.
\end{example}

\begin{example}[Dependent linkage classes]
Consider the weakly reversible network with species $X,A,B$ and reactions
\[
X\rightleftharpoons A\rightleftharpoons B,\qquad X+B\rightleftharpoons 2A .
\]
It has five complexes and two linkage classes, $K_1=\{X,A,B\}$ and $K_2=\{X+B,2A\}$. The reaction vectors span the plane $S=\{v: v_1+v_2+v_3=0\}$, so $s=2$, while $\dim S_{K_1}=2$ and $\dim S_{K_2}=1$. Since $2+1>2$, the linkage classes are stoichiometrically dependent, and the deficiency is $5-2-2=1$. Here $S^\perp=\spn(1,1,1)$, the projection of $e_X$ is $\psi=\tfrac13(1,1,1)$, and $\psi_{K_1}=\tfrac13$, $\psi_{K_2}=\tfrac23$. Theorem~\ref{thm:main} predicts $(z_t)_X\asymp e^{t/3}$.

We checked this numerically. The rate constants were $(1,2)$ for $X\rightleftharpoons A$, $(0.7,1.3)$ for $A\rightleftharpoons B$ and $(3,0.4)$ for $X+B\rightleftharpoons 2A$ (forward, backward), and $b=(1,2,0.5)$. We applied $40$ random rescalings, with an independent factor in $[e^{-6},e^{6}]$ for each of the three reversible pairs. This is more general than Proposition~\ref{prop:uniform} requires, and the proof applies verbatim to any partition of the reactions into strongly connected pieces. All computed steady states in $P_b$ lay in $[0.076,\,2.28]^3$. For $t=-6,-3,0,3,6$, the points $z_t$ were steady states of the original system with $x_X=0.268,\,0.728,\,1.98,\,5.38,\,14.6$, increasing by a factor $e$ each time $t$ increases by $3$.
\end{example}

\section{Discussion}\label{sec:discussion}

Corollary~\ref{cor:acr} says that nontrivial ACR, in any mass-action system and for any rate constants, requires a reaction network that is not weakly reversible, that is, a reaction whose source cannot be recovered from its target. In biochemical models such steps represent processes held far from equilibrium, typically by ATP hydrolysis. This is consistent with the observation that the Shinar--Feinberg mechanism relies on nonterminal complexes.

The result fits with stochastic considerations. Anderson, Enciso and Johnston showed that the Shinar--Feinberg ACR networks, modelled as continuous-time Markov chains, undergo extinction with probability one \cite{AEJ14}. For weakly reversible networks, on the other hand, Anderson and Kim conjectured that the stochastic mass-action model is positive recurrent on every closed communicating class \cite{AK18}. A proof of this conjecture is given in a companion paper [in preparation]. Together with Corollary~\ref{cor:acr}, it implies that the structural property guaranteeing stochastic persistence for all rate constants excludes nontrivial ACR for all rate constants.

\subsection*{Formal verification}
Proposition~\ref{prop:uniform}, Lemma~\ref{lem:shift}, Theorem~\ref{thm:main} and Corollary~\ref{cor:acr} have been formalised in Lean~4 with Mathlib, for finite mass-action networks over $\R$ with an arbitrary labelling of linkage classes. The existence theorem of Boros (Theorem~\ref{thm:boros}) enters as an explicit hypothesis and is not formalised. The file \texttt{WACR.lean} is included with the source of this paper, is available at \url{https://github.com/kuhaizhongdebailanren/AIM272_lean}, and is archived at \url{https://doi.org/10.5281/zenodo.24503870} (release v1.0.0, folder \texttt{WACR}) \cite{PDZlean}. It compiles with Lean~4.24.0 and the corresponding Mathlib, and the three main theorems depend only on the standard axioms \texttt{propext}, \texttt{Classical.choice} and \texttt{Quot.sound}.

\subsection*{Statement on the use of AI tools}
In accordance with arXiv's policy on generative AI tools, we report that large language models developed by OpenAI (GPT-6.2 Sol and GPT-6.2 Astra) were used in developing and checking the arguments of this paper. All mathematical content generated with the help of these tools was reviewed, revised and checked by the authors. In addition, the results listed in the paragraph on formal verification above have been machine-checked in Lean~4, with Theorem~\ref{thm:boros} as an explicit hypothesis; the code is archived at \url{https://doi.org/10.5281/zenodo.24503870}. AI tools are not authors of this paper. The authors take full responsibility for its entire content, including any errors.


\begin{thebibliography}{99}

\bibitem{AEJ14}
D.~F. Anderson, G.~A. Enciso, and M.~D. Johnston,
\emph{Stochastic analysis of biochemical reaction networks with absolute concentration robustness},
J. R. Soc. Interface \textbf{11} (2014), no.~93, 20130943. doi:10.1098/rsif.2013.0943.

\bibitem{AK18}
D.~F. Anderson and J.~Kim,
\emph{Some network conditions for positive recurrence of stochastically modeled reaction networks},
SIAM J. Appl. Math. \textbf{78} (2018), no.~5, 2692--2713.

\bibitem{BG03}
E.~Batchelor and M.~Goulian,
\emph{Robustness and the cycle of phosphorylation and dephosphorylation in a two-component regulatory system},
Proc. Natl. Acad. Sci. USA \textbf{100} (2003), no.~2, 691--696.

\bibitem{Bor19}
B.~Boros,
\emph{Existence of positive steady states for weakly reversible mass-action systems},
SIAM J. Math. Anal. \textbf{51} (2019), no.~1, 435--449. doi:10.1137/17M115534X.

\bibitem{HJ72}
F.~Horn and R.~Jackson,
\emph{General mass action kinetics},
Arch. Ration. Mech. Anal. \textbf{47} (1972), 81--116. doi:10.1007/BF00251225.

\bibitem{JC22}
B.~Joshi and G.~Craciun,
\emph{Foundations of static and dynamic absolute concentration robustness},
J. Math. Biol. \textbf{85} (2022), Paper No.~53. doi:10.1007/s00285-022-01823-2.

\bibitem{KNS24}
N.~Kaihnsa, T.~Nguyen, and A.~Shiu,
\emph{Absolute concentration robustness and multistationarity in reaction networks: conditions for coexistence},
European J. Appl. Math. \textbf{35} (2024), no.~4, 566--600. doi:10.1017/S0956792523000335.

\bibitem{MST22}
N.~Meshkat, A.~Shiu, and A.~Torres,
\emph{Absolute concentration robustness in networks with low-dimensional stoichiometric subspace},
Vietnam J. Math. \textbf{50} (2022), no.~3, 623--651. doi:10.1007/s10013-021-00524-5.

\bibitem{PDZlean}
S.~Peng, C.~Deng, Y.~Zou, and R.~Li,
\emph{AIM272\_lean: Lean~4 formalization},
Zenodo, 2026. doi:10.5281/zenodo.24503870.

\bibitem{SF10}
G.~Shinar and M.~Feinberg,
\emph{Structural sources of robustness in biochemical reaction networks},
Science \textbf{327} (2010), 1389--1391. doi:10.1126/science.1183372.

\end{thebibliography}
\end{document}